\documentclass[%
 aip,
 sd,%
 amsmath,amssymb,
 preprint,%
]{revtex4-1}

\usepackage{graphicx}
\usepackage{dcolumn}
\usepackage{subfigure}
\usepackage{bm}
\usepackage{color}

\begin{document}

\preprint{AIP/123-QED}

\title[Refining and classifying FTLE ridges]{Refining and classifying finite-time Lyapunov exponent ridges}

\author{M. R. Allshouse}
 \email{mallshouse@chaos.utexas.edu.}
 \altaffiliation[Now at ]{Physics Department, University of Texas - Austin.}
\author{T. Peacock}%
 \email{tomp@mit.edu.}
\affiliation{ 
Mechanical Engineering Department, Massachusetts Institute of Technology
}%

\date{\today}

\begin{abstract}
While more rigorous and sophisticated methods for identifying Lagrangian based coherent structures exist, the finite-time Lyapunov exponent (FTLE) field remains a straightforward and popular method for gaining some insight into transport by complex, time-dependent two-dimensional flows. In light of its enduring appeal, and in support of good practice, we begin by investigating the effects of discretization and noise on two numerical approaches for calculating the FTLE field. A practical method to extract and refine FTLE ridges in two-dimensional flows, which builds on previous methods, is then presented.  Seeking to better ascertain the role of an FTLE ridge in flow transport, we adapt an existing classification scheme and provide a thorough treatment of the challenges of classifying the types of deformation represented by an FTLE ridge. As a practical demonstration, the methods are applied to an ocean surface velocity field data set generated by a numerical model.  

%
\end{abstract}

\maketitle

\begin{quotation}
The transport of material by spatio-temporally complex flow fields has widespread application to geophysical and industrial processes. In recent years, advances in methods to detect Lagrangian Coherent Structures (LCS), these being key material lines within the flow field that organize flow transport, have enabled significant breakthroughs and profound new insight. One of the earliest developed LCS detection methods employs the finite-time Lyapunov exponent field (FTLE), and while newer and more rigorous methods exist, the simplicity of the FTLE approach, both in its implementation and its identification of  active regions of a flow field, means that it is still widely used. As such, the thrust of this work is to highlight good practices in calculating the FTLE field, to provide a practical means to identify the most active material lines of a flow field that lie along maximal FTLE ridges, and to classify the types of deformation associated with an FTLE ridge. The methods we refine and advocate, which build on previous established results, are tested on an analytical model and a numerical model data set of an ocean surface flow.
 
\end{quotation}

\section{\label{sec:intro}Introduction}

There has been substantial recent effort to develop, establish and advance rigorous, objective methodologies and associated data analysis tools for understanding how advective transport is organized in complex, time dependent fluid flows~\cite{Peacock2013, Samelson2013, Haller2015}. One of the first tools developed for this purpose was the Finite Time Lyapunov Exponent (FTLE), which reveals the regions of the flow field that undergo the greatest stretching for the time window considered~\cite{Haller2001, Shadden2005}. More sophisticated and rigorous approaches have subsequently been developed; in particular, geodesic methods identify structures such as strainlines and shearlines, whose properties are very well defined and for which there are numerically stable tools~\cite{Haller2012}. The progress has been substantial, particularly for two-dimensional flows that are relevant to such important scenarios as the evolution of an oil spill on the ocean surface~\cite{Olascoaga2012, Allshouse2015}.

Although more advanced methods now exist, obtaining the FTLE field and visualizing the FTLE ridges remains a popular and insightful means for investigating the organization of transport in complex flows.  A healthy degree of caution and understanding are needed when interpreting the results.  For example, it should be recognized that since FTLE ridges are simply representations of the time history of what happens to material elements, it may be that an entire FTLE ridge exists only because the material elements of which it is comprised are transported past an important local Lagrangian flow feature (e.g. a hyperbolic core)~\cite{Olascoaga2012, Allshouse2015}. It is unarguably the case, however, that FTLE ridges represent the most kinematically active (i.e. most local stretching) material lines of the flow field, and this may be the most important consideration for a study.  Furthermore, strainlines and shearlines, these being key material lines identified by the geodesic approach~\cite{Haller2012}, have no obligation to align with FTLE ridges in general, compressible two-dimensional flows, and so one cannot assume that an FTLE ridge will be so marked (although it is always worth checking).  

While FTLE analysis has been, and continues to be, widely utilized~\cite{Shadden2005, Mathur2007}, there is scope for improvement.  For example, as with any numerical method, it is important to check convergence of the FTLE values with system parameters, such as  the velocity field resolution or the cluster size used for finite-difference approximation; while somewhat rudimentary, this topic has yet to be systematically addressed in the FTLE literature.  It also remains a practical challenge to robustly extract ridges of the FTLE field,  which provides a simpler presentation of the results and enables material line advection, thereby revealing  how FTLE ridges evolve and shape transport.  Finally, an outstanding issue is that FLTE ridges do not reveal what type of local deformation they represent (i.e. normally hyperbolic repulsion, Lagrangian shear or tangential stretching) and there has been only one previous study~\cite{Tang2011} that attempted to take this further step; the ability to do so would provide more insight into the nature and significance of an FTLE ridge.

In this paper, recognizing the continued use by many of FTLE analysis particularly when studying geophysical flows, we revisit the topic and seek to address the aforementioned outstanding issues.  In section~\ref{sec:models}, we demonstrate good practices in checking convergence of the FTLE field, as well as outlining an alternative method for calculating FTLE fields that is highly accurate for analytic data sets.  Next, section~\ref{sec:ridges} presents a practical extraction and refinement scheme for determining FTLE ridges with sufficient accuracy that, despite being the most sensitive features in the flow field, they can be faithfully advected.  Section~\ref{sec:classification} then attempts to classify FTLE ridges according the type of  local deformations they represent, identifying the inherent challenges and resolving what classifications are reasonably achievable and under what circumstances.  Having addressed these issues using an analytic model as a test case, in section~\ref{sec:ningaloo} we proceed to apply our findings to a real-world case study of an ocean surface flow, using a data set generated by a high quality numerical model.  Finally, in section~\ref{sec:conclusions}, we present our conclusions. Some technical details are presented in the appendices.

\section{FTLE Calculations}
\label{sec:models}

The principal task underlying the determination of the FTLE field is the calculation of the right Cauchy-Green strain tensor (henceforth referred to as the CG tensor) for the physical domain and time window of interest. In this section, we discuss the standard  finite-difference methods and a less well utilized approach to calculate the CG tensor, investigate the impacts of discretization and noise, and discuss good practices.  While the focus of the present study is FTLE ridges, it should be noted that the following calculations are also necessary for strainline and shearline based methods.

\subsection{Methods}

The first step in calculating the CG tensor is to calculate the flow map, which maps a fluid element from its initial position $\bm{x}_0=(a_1,a_2)$ at time $t_0$ to its final position $\bm{x}=(x_1,x_2)$ at time $t$, and is represented as $\bm{F}^{t}_{t_0}(\bm{x}_0) = \bm{x}(t;\bm{x}_0, t_0)$. To elucidate the notation, the initial condition, also the Lagrangian coordinates, $(a_1,a_2)$ are used to differentiate from the Eularian coordinates $(x_1, x_2)$.  Given a velocity field $\bm{u}(\bm{x},t)$, solutions to the flow map are determined by solving the system of Ordinary Differential Equations (ODEs):
\begin{equation}
\label{eq:vel_eom}
\frac{d \bm{x}}{d t} = \bm{u}(\bm{x},t),
\end{equation}
for an initial position of $\bm{x}_0$ at time $t_0$.  We use the Matlab solver {\tt ode45}, which is a variable step Runge-Kutta 4$^{th}$ and 5$^{th}$ order ODE solver, for all calculations in this paper. If the velocity field is analytic, then velocities can be calculated directly.  If it is a discrete velocity field data set, however, velocities at off-grid locations required by the ODE solver must be obtained by interpolation, the impacts of which are investigated later in this section. 

 In order to map a set of initial points $\{\bm{x}^{(1)}_0,\bm{x}^{(2)}_0,...,\bm{x}^{(n)}_0\}$, the equations of motion \eqref{eq:vel_eom} must be solved for each initial condition.  This process is accelerated by simultaneously solving the series of equations:
\begin{align}
\label{eq:group}
\dot{\bm{x}}^{(1)}(\bm{x}^{(1)},t) &= \bm{u}(\bm{x}^{(1)},t), \nonumber \\
\dot{\bm{x}}^{(2)}(\bm{x}^{(2)},t) &= \bm{u}(\bm{x}^{(2)},t), \nonumber \\
\cdots &= \cdots, \nonumber \\
\dot{\bm{x}}^{(n)}(\bm{x}^{(n)},t) &= \bm{u}(\bm{x}^{(n)},t).
\end{align}
Simultaneous advection takes advantage of the built in vectorization of Matlab and imposes the same variable time step sequence and numerical tolerances to all trajectories. For sensitive systems, failure to advect trajectories simultaneously may result in erroneous results, particularly if the relative tolerances are large.  As good practice, the absolute and relative tolerance for the solver should be systematically reduced so that convergence of the trajectories is achieved at the associated but necessary cost of increased computational run time.

For a given initial condition $\bm{x}_0$, the advected final position of a nearby point, $\bm{x}_0+\bm{\epsilon}$, can be approximated by:
\begin{equation}
\bm{F}^t_{t_0}(\bm{x}_0+\bm{\epsilon}) = \bm{F}^t_{t_0}(\bm{x}_0) + \bm{\nabla F}^t_{t_0}(\bm{x}_0) \bm{\epsilon} + \textit{O}(|\bm{\epsilon}|^2), 
\end{equation}
where the second term in this expansion contains the flow map gradient:
\begin{equation}
\bm{\nabla F}^t_{t_0}(\bm{x}_0) = \partial x_i / \partial a_j\Big|_{\bm{x}_0} \mbox{ for } i,j = 1,2,
\end{equation}
Finite-difference methods have been the primary means for calculating the flow map gradient, with nearest neighbors on a regular grid initially being used for the derivative calculations~\citep{Lekien2005}.  Uniform-grid based methods pose an issue when trying to improve the accuracy of the finite-difference scheme, however, due to the rapidly increasing computational demand as the grid spacing is reduced.  Unstructured mesh methods have been utilized, with dynamic mesh refinement used to improve the resolution of flow map calculations in high FTLE regions~\citep{Lekien2010}.  While this improves accuracy, the need for mesh refinement complicates the calculation, increasing the computational demand.  The use of clusters~\cite{Farazmand2012}, which adds four off grid points of adjustable distance to calculate the finite-difference term for each on-grid location, has become widely used.  This method allows for Cartesian grids with run times independent of the desired accuracy, at the not unreasonable expense of a five-fold increase in computational cost.   Arbitrarily fine accuracy cannot be achieved as the cluster size is reduced, of course, as the difference in advected position eventually becomes the same order as the numerical resolution, resulting in errors.  Underlying all these finite-difference techniques is the possibility that neighboring points may not remain nearby throughout advection. In such a situation, the finite-difference approximation of a local deformation may itself not be valid.  To compensate, one can renormalize the cluster around the central point throughout advection preventing nearby trajectories from becoming non-local~\cite{Nese1989}; this, of course, increases computational complexity and demands. 

An alternative, but relatively untested, approach for calculating the terms of the flow map gradient is to use the advected-gradient method, which solves a system of ODEs that yields the terms of both $\bm{F}^t_{t_0}(\bm{x})$ and $\bm{\nabla F}^t_{t_0}(\bm{x}_0)$ directly. This is achieved by solving~\eqref{eq:vel_eom} and the following system of equations, obtained via the chain rule, as follows:
\begin{equation}
\frac{d}{dt}\frac{\partial x_i}{\partial a_j} = \frac{\partial}{\partial a_j}\frac{d x_i}{dt} = \frac{\partial x_k}{\partial a_j}\frac{\partial u_i}{\partial x_k} \mbox{ for } i,j=1,2,
\end{equation}
where we have used standard index notation convention for summing repeated indices; the initial conditions for these equations are: 
\begin{equation}
\bm{x}(t_0) = \bm{x}_0, \hspace{1cm} \frac{\partial x_1}{\partial a_1} = \frac{\partial x_2}{\partial a_2} = 1, \hspace{1cm} \frac{\partial x_1}{\partial a_2} = \frac{\partial x_2}{\partial a_1} = 0.
\end{equation}
An advantage of the advected-gradient approach is that for an analytic system, where the velocity field and its gradients are completely known, this method will be accurate.  The six coupled equations typically take longer to solve than the ten equations needed for a cluster based finite-difference approach (two equations for each of the five particles in a cluster); this is due to the large and rapidly changing values of the flow-map-gradient terms requiring more computation time to satisfy equivalent numerical tolerances.  It also remains to be determined to what extent the advected-gradient method is compromised by discretization and noise.

Once the flow map gradient field has been determined, the FTLE field is readily extracted from the CG tensor field:
\begin{equation}
\bm{C}=[\bm{\nabla F}^t_{t_0}]^{*}[\bm{\nabla F}^t_{t_0}],
\end{equation}
where $^*$ corresponds to the transpose operator. By definition, the FTLE is:
\begin{equation}
\Phi = \frac{1}{2(t-t_0)} \log \lambda_2,
\end{equation}
where $\lambda_2$ is the largest eigenvalue of $\bm{C}$.

\subsection{Results}
\label{sec:ftle_calc}

To test the FTLE methods, we use an autonomous analytic model that has analytic solutions for all the quantities we are concerned with, so that in all cases there is a `true' result against which to compare the results of our various numerical calculations. Here we outline the model, more detail on which is given in appendix~\ref{app:model}. The autonomous analytical system is based on the nonlinear system of equations:
\begin{align}
\label{eq:base_equations}
\dot{X_1} &= -X_1{}^3 + X_1, \nonumber \\
\dot{X_2} &= X_2{}^3 - X_2,
\end{align}
within the domain $D=\{ X_1 \in [-1,1], X_2 \in [-1,1] \}$. To add spatial complexity to the flow field we introduce a coordinate transformation:
\begin{align}
\label{eq:transform_cord}
x_1 &= X_1 \cos(r) - X_2\sin(r), \nonumber \\
x_2 &= X_2 \cos(r) + X_1\sin(r), 
\end{align}
where $r = \sqrt{x_1{}^2+x_2{}^2} = \sqrt{X_1{}^2+X_2{}^2}$. The equations of motion for the transformed coordinates are calculated by taking the time derivative of equation~\eqref{eq:transform_cord}:
\begin{align}
\label{eq:swirl_eom}
\dot{x_1} &= \left( x_1 - x_1^3 - (\frac{1}{r}-r)(x_1^2-x_2^2)x_2 \right) \cos(2r) \nonumber\\
& \hspace{1cm} +\left(x_2-\frac{1}{2}x_2(3x_1^2+x_2^2) - (\frac{1}{r}-r)2x_1x_2^2\right) \sin(2r), \nonumber\\
\dot{x_2} &= \left(-x_2 + x_2^3 + (\frac{1}{r}-r)(x_1^2-x_2^2)x_1 \right) \cos(2r) \nonumber \\
& \hspace{1cm} + \left(x_1-\frac{1}{2}x_1(x_1^2+3x_2^2) - (\frac{1}{r}-r)2x_1^2x_2\right) \sin(2r).
\end{align}
Due to the transformation, the repelling invariant manifolds, $X_2=\pm1$, now correspond to the left and right boundaries. The origin remains a stationary point in phase space and is still a hyperbolic core~\cite{Olascoaga2012}.  Neither the original system nor the transformed system satisfy continuity, which is intentionally so in order to permit a greater variety of local deformations than possible for an incompressible system.  While the level of compressibility in the model system is strong, in the ocean, for example, upwellings and downwellings can cause the surface velocity field (i.e. not accounting for the vertical velocity) to be highly compressible.

\begin{figure}
\centering
		\includegraphics[width=6.5in]{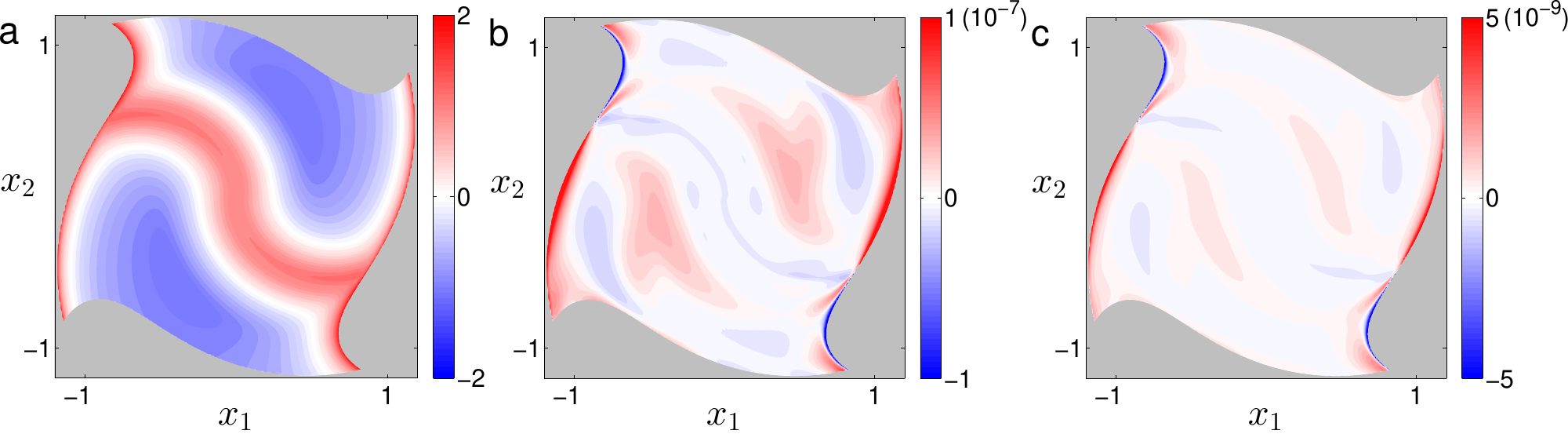}
		\caption{(a) Analytic FTLE field for the autonomous system with $t-t_0=2$.  (b) Difference between the analytic and finite-difference FTLE fields.  (c) Difference between the analytic and advected-gradient FTLE fields. }
	\label{fig:swirl_ftle}
\end{figure}

The analytic FTLE field for a time window $t-t_0=2$ is presented in figure~\ref{fig:swirl_ftle}(a).  The largest FTLE values are along the left and right boundaries and through the center of the domain.  While the largest values of the FTLE field exceed 2, there is a large portion of the domain where the FTLE value is negative, highlighting the fact that the system does not satisfy continuity. For comparison, the error in the FTLE fields calculated via the finite-difference and advected-gradient approaches are presented in figures \ref{fig:swirl_ftle}(b) and (c), respectively.  The finite-difference calculation was most accurate when using a cluster size $\delta a=10^{-6}$ and a relative tolerance of $10^{-7}$; these parameters yielded a spatially averaged error on the order of $10^{-8}$, with the largest errors occurring along the curvilinear left and right boundaries.  While the errors in the FTLE values along the boundary are both positive and negative, there is a persistent underestimation of the FTLE field across the center of the domain.  The advected-gradient result has a similar structure to the finite-difference error field; the former is on average two orders of magnitude smaller, however, with a spatially-average error of order $10^{-10}$.  For both calculations, identical relative tolerances were used, therefore the additional error is solely due to the finite-difference approximation.

Typically, for FTLE calculations the analytic velocity field is not known and instead a velocity field data set discretized in space and time is available.  It has been established that FTLE ridges are robust to isolated errors in the velocity field either spatially or temporally given that the error is small~\cite{Haller2001}, but persistent errors resulting from discretization or noise may not satisfy this condition.  When analyzing ocean models, spatial resolution is fundamental to determining the scale of FTLE structures that are reliably detected.  A number of studies have investigated the effects that ocean model resolution has on the resulting FTLE or finite-size Lyapunov exponent fields~\cite{BeronVera2010, Hernandez2011, Poje2010}.  One particular study systematically subsampled turbulent ocean models to study the impact of spatial and temporal on the FTLE calculation~\cite{Keating2011}.  Another study has investigated the impact of noise and discritization on an FTLE field calculated via finite-difference on a uniform field~\cite{Olcay2010}; here we focus on the relative impacts of both on the different FTLE calculation methods.  Our goal is to study the impact of varying spatial resolution and the presence of noise on the finite-difference based FTLE calculation relative to the analytically calculated values.  

\begin{figure}
	\centering
		\includegraphics[height=2in]{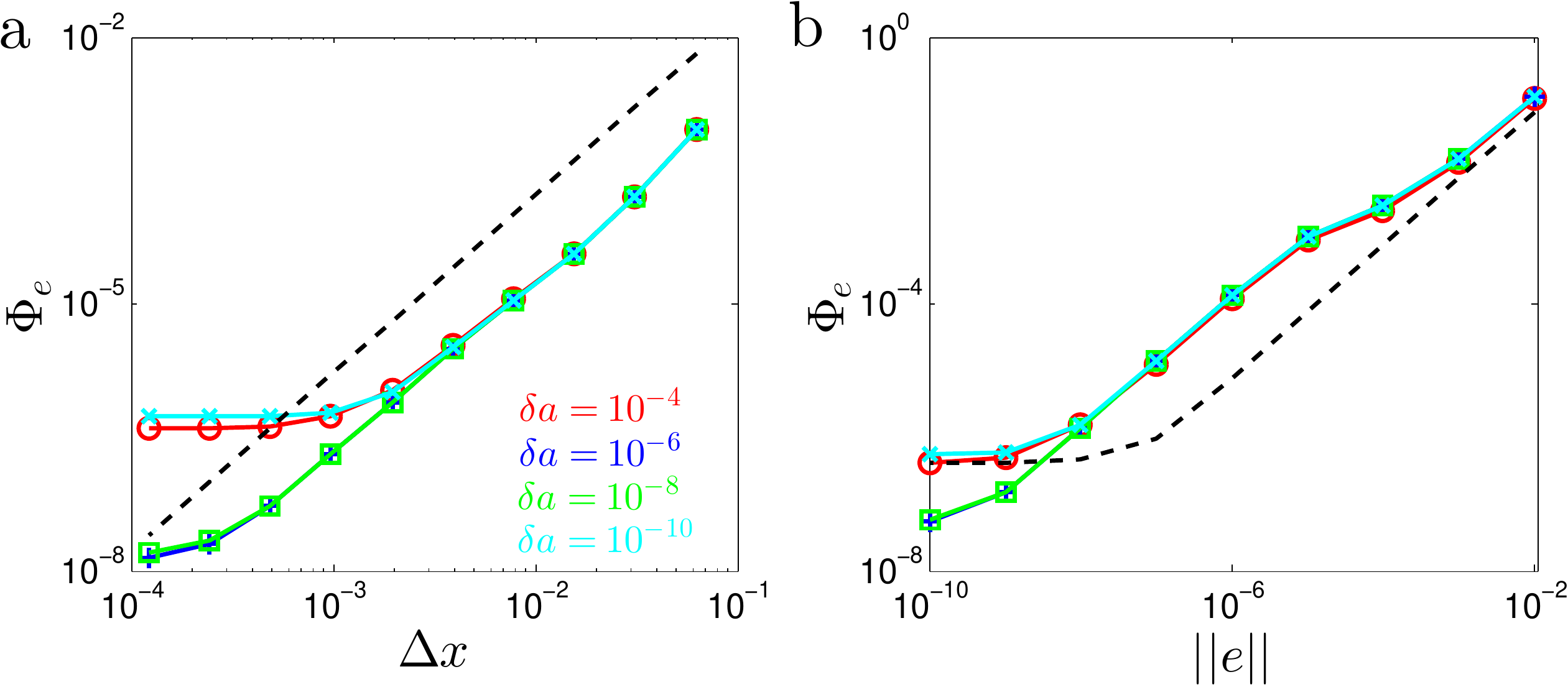}
		\caption{(a) Spatially averaged FTLE error for regions where $\Phi\geq 1$ for the advected-gradient (dashed) and finite-difference methods (solid) calculated from a discrete velocity field with varying resolution, $\Delta x$; the cluster spacing, $\delta\/a$, is also varied.  All calculations used a relative tolerance of $10^{-7}$ for the ODE solver.  (b)  Error calculated from a discrete velocity field with a resolution $\Delta x=2^{-11}$ with uniformly distributed noise of maximum magnitude $||e||$. }
	\label{fig:comp_study}
\end{figure}

To investigate the effects of discretization, we considered a uniform grid of data from our analytical model for the range of grid spacings $\Delta x = 2^{-4}$ to $2^{-13}$, utilizing Matlab's {\tt griddedInterpolants} to accelerate the interpolation. Because Matlab's interpolant is faster when performing simultaneous as opposed to sequential interpolations, the simultaneous advection scheme given by equation \eqref{eq:group} again has significant advantages over the advection of individual trajectories. To quantify the degree of error in regions of interest, we calculated the spatially averaged error of the FTLE field, $\Phi_e$, in regions where $\Phi\geq1$.  This measure is presented as a function of the velocity field resolution for both the finite-difference and advected-gradient methods in figure~\ref{fig:comp_study}(a); furthermore, for the finite-difference method the results are obtained for cluster sizes $\delta\/a=10^{-4}, 10^{-6}, 10^{-8}, $ and $10^{-10}$.  Both methods show order $0.1\%$ errors for the coarsest velocity field resolution.  As $\Delta\/x$ is reduced, the finite-difference results converge to the numerical errors obtained when working with the analytic velocity field.  The error of the advected-gradient method, however, is typically two orders of magnitude larger than those of the finite-difference method.  This is attributed to the sensitivity of the flow map gradient terms to errors in the velocity gradient field, and reveals that the finite-difference method performs better for discretized data. 

In addition to discrete data sets, another concern is that velocity field data sets may be noisy. To systematically study the effect of noise, a data set with spatial resolution $\Delta x=2^{-11}$ was used (so that noise was the dominant factor over discretization) and at each grid point uniformly distributed noise of magnitude $||e||$ was added to the velocity field, which is of order $10^{-1}$.  The results for $\Phi_e$ for the two methods are presented in figure~\ref{fig:comp_study}(b).  For large values of $||e||$, the advected-gradient and finite-difference methods are similarly impacted by noise, with the former method somewhat outperforming the latter as noise levels are reduced. 

Based on these studies, the introduction of a discretized field, with or without noise, effectively eradicates the greater accuracy of the advected-gradient method, and thus the cluster-based finite-difference method should be employed.  While this study has been performed for a simple analytic system, we expect these results to carry over to more spatially-complex, time-dependent discretized velocity fields, as we performed similar studies for the double-gyre flow that generated similar results.  

\section{FTLE Ridges}
\label{sec:ridges}

The standout features of an FTLE field are the maximal ridges, which correspond to regions that undergo maximum Lagrangian separation of neighboring fluid elements. As has been pointed out, FTLE ridges do not necessarily coincide with variational\cite{Haller2011} and geodesic~\cite{Haller2012} definitions of LCS, but they may nevertheless be considered a type of LCS, in the sense that they identify the most kinematically active material lines in a flow field for a given time window. It is desirable to be able to accurately identify the material lines that reside along FTLE ridges in order to form a simplified picture of the flow transport, to enable these structures to be advected over the time window, and, as we shall see, to enable classification of the associated types of deformation. In so doing, it should be appreciated that being in such active regions of the flow field, these material lines will likely be greatly deformed as they are advected and may not represent a barrier to transport, per se.

\subsection{Methods}
\label{sec:ridge_def}

To develop the foundation for detecting and classifying FTLE ridges, it is first necessary to define a ridge. For two-dimensional systems, ridges of the FTLE field are one-dimensional lines that mark generalized local maxima of the field.  A ``height ridge'' is a line made up of a set of points that are locally maximum in the direction of the greatest curvature of the FTLE field  projection~\cite{Eberly1996}. While this definition lends itself to a point by point evaluation, it does not necessarily create connected lines forming ridges, and attempts to amend the height ridge definition\cite{Shadden2005} resulted in an over-constrained system~\cite{Schindler2012}.  An alternative to the height ridge, the ``watershed ridge'' divides the system into disjoint regions based on global stationary points of the field~\cite{Eberly1996}.  The ridges are identified as slope lines, these being trajectories of the FTLE gradient field that connect the saddle points in the system to the local maxima. One variation of the watershed ridge removed the requirement that it start at stationary points of the FTLE field but imposes a requirement that the slope line be a normally attracting invariant manifold~\cite{Karrasch2013}.  Both the height ridge and watershed ridges definitions are practically challenging to implement in the presence of noisy data. In the case of the height ridge, the challenge arises when trying to calculate the Hessian matrix of the FTLE field.  And while the watershed definition does naturally produce connected lines, noise may produce many stationary points in the FTLE field and significantly contaminate the gradient field.

 Building on the definitions of height ridges and watershed ridges, we form a numerically tractable ridge definition.  Principally, a ridge should be a set of connected points that are generalized local maxima of the FTLE field, so points on the ridge will be a local maximum FTLE value in the direction normal to the ridge.  Compared to the height ridge definition, this definition loosens the requirement that the ridge normal be aligned with the smallest eigenvector of the Hessian of the FTLE field, thus removing the need to calculate second derivatives of the numerically evaluated FTLE field.  To form an actual line, ridges will be everywhere tangent to the gradient of the FTLE field, like the slope lines that form watershed ridges, but there is no requirement that the starting points of the ridges be stationary points of the FTLE field, which removes the need to accurately identify the numerically evaluated stationary points.  With these criteria in hand, we define a {\it normal-maximum ridge} as the line $\gamma$ where $\forall \bm{x} \in \gamma$:
\begin{subequations}
\begin{align}
\label{eq:ridge_crit1}
\bm{e}(\bm{x}) &= \frac{\bm{\nabla}\Phi}{|\bm{\nabla}\Phi|}, \\
\label{eq:ridge_crit2}
\bm{n}(\bm{x})^*\bm{H}_{\Phi}(\bm{x})\bm{n}(\bm{x}) &< 0,
\end{align}
\end{subequations}
where $\bm{e}(\bm{x})$ and $\bm{n}(\bm{x})$ are the tangent and normal to $\gamma$ at $\bm{x}$, respectively.  The first condition ensures that the ridges are slope lines, and the second condition ensures points are a local maximum in the direction normal to the ridge.  In practice, the Hessian is not actually calculated, and instead  FTLE values of points nearby and along a normal to the ridge are used to reasonably check the local maximum criteria.  

\begin{figure}
	\centering
		\includegraphics[height=3in]{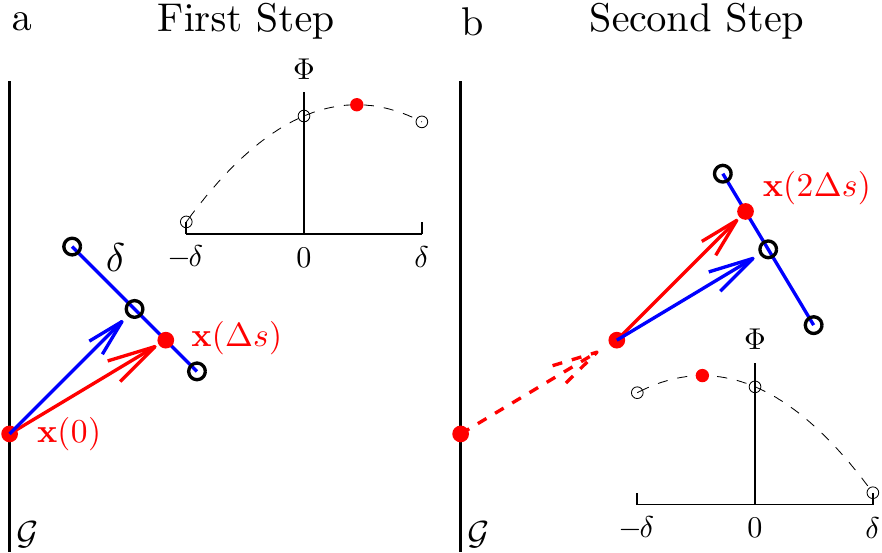}
		\caption{(a) Initial step of the ridge tracking.  The seed point (red) is a local max along the grid $\mathcal{G}$.  An initial step direction is set by the blue arrow. The set of test points are black circles with the FTLE values as shown in the inset.  The corresponding maximum is shown as a red dot in both the main and inset plots.  (b) Second step in the ridge tracking process.  The direction is set by the vector from the seed point and the initial step (red dashed arrow).}
	\label{fig:ridge_tracking}
\end{figure}

There already exists a numerical algorithm that with minor modifications can be used to locate normal-maximum ridges~\cite{Shadden2005, Lipinski2010}. The basic ridge tracking algorithm, as illustrated in figure~\ref{fig:ridge_tracking}, is as follows:
\begin{enumerate}
\item Locate seed points, $\bm{x}(s=0)$, which are local FTLE maxima in the initial step direction, $\bm{\nabla}\Phi$.  To accelerate the search, points are selected from a grid, $\mathcal{G}$, made up of vertical and horizontal lines that divide the system.  Start two trajectories from each seed point with initial step directions opposite of each other
\item Step in this direction a distance $\Delta s$, and calculate the FTLE values at the new position as well as two other points, one either side in a direction normal to the step direction taken.
\item Fit a parabola to the FTLE values calculated and add the corresponding location of the maximum, $\bm{x}(s+\Delta s)$ to the ridge.
\item Use the previous step position, $\bm{x}(s)$, and the current position to calculate the tangent vector.
\item Repeat steps 2-4 until a stop condition is met, such as the next step (i) leaves the domain, (ii) hits the start or end of another ridge, (iii) fails to find a maximum or (iv) has an FTLE value below a user defined threshold.
\end{enumerate}
The fundamental difference between this approach and the previous version is how the FTLE field is calculated~\cite{Lipinski2010}.  To reduce the number of calculations, the previous version sequentially calculates the FTLE field at test points, ignoring regions of the field where there are no FTLE ridges.  Because the simultaneous calculation of the FTLE field is more efficient than sequential calculation, however, we calculate the FTLE field for the entire domain and use interpolation to calculate the FTLE values for the above algorithm. This proves sufficient to find a reasonable first approximation of the ridge, which can then be refined.

For systems where the FTLE field has sharp ridges, and thus there are large magnitudes of the second-derivative in the direction normal to the ridge, it is more appropriate to select an initial step in the direction normal to the gradient.  This is due to large gradients normal to the ridge for positions not precisely on the ridge~\cite{Haller2001}.  For subsequent steps, the approach does not rely on the gradient calculation, so the method accurately tracks the FTLE ridges.  Regardless of ridge sharpness, errors in the ridge tracking step are particularly large in the initial steps due to not being directly on the ridge.  This error is reduced by the refinement approach discussed next.

\begin{figure}
\centering
		\includegraphics[height=2in]{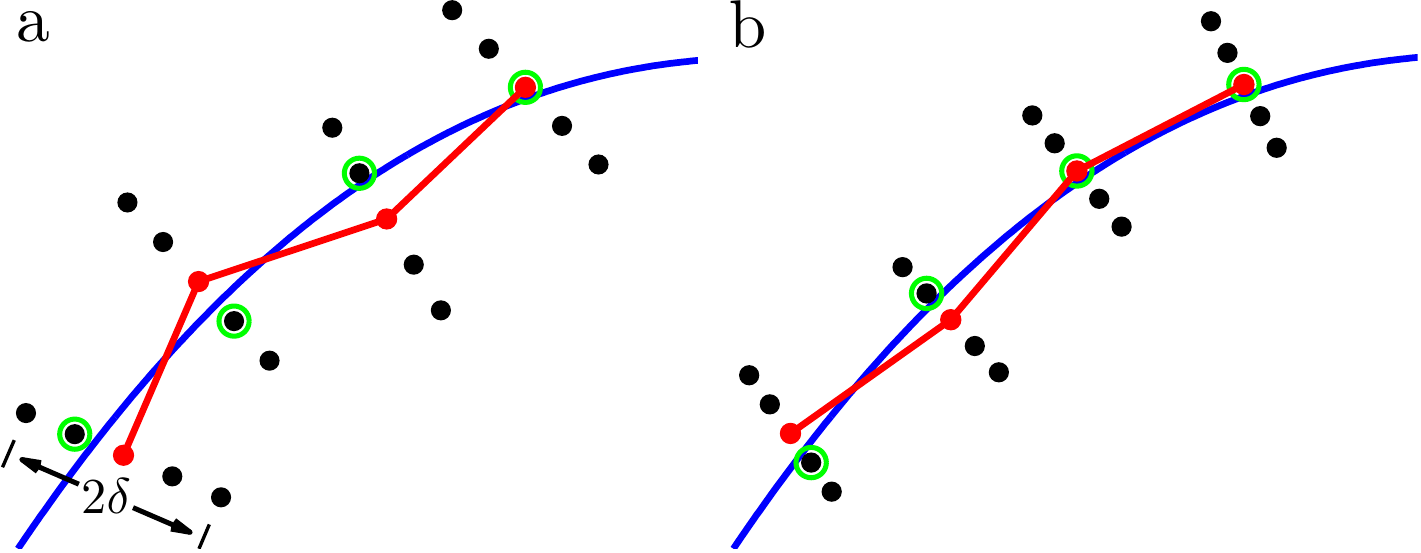}
		\caption{(a) Initial ridge detected via ridge tracking (red) with the points normal to the ridge (black).  The points closest to the true ridge (blue) have the highest FTLE values and are circled (green).  (b) The updated position of the estimated ridge using the circled points in (a) is in red.  Normal points spanning a smaller search range (black) are again calculated.  The new maximum position is circled in green.}
	\label{fig:ridge_refine}
\end{figure}

Having used the above recipe to initially locate a ridge, a refinement scheme is then employed.  As depicted in figure~\ref{fig:ridge_refine}(a), the refinement scheme takes an initial ridge (red), calculates the normal at all points, and then places a number of test points at incremental distances $\delta$ either side of the ridge, along the normal direction (black points).  FTLE values for all these points normal to the ridge are calculated (not interpolated) and the point with the maximum FTLE value at each cross ridge location is taken as the revised position of the ridge (green circles).  For the updated position of the ridge, the normal along the ridge is recalculated, a smaller search range either side of the revised ridge is considered, the FTLE values for points in this range are calculated, and the ridge position is updated, as illustrated in figure~\ref{fig:ridge_refine}(b).  This process is repeated iteratively for progressively smaller values of $\delta$ until the accepted degree of accuracy is reached.  

One of the principal benefits of this method over previous methods is that we find that an initial calculation of the entire FTLE field is significantly faster than running the FTLE calculation of a small set of points hundreds of times; for even moderate resolution, the interpolated based ridge tracking still produces an acceptable first estimate for the ridges of the FTLE field, which is the initial guess needed to proceed to the ridge refinement scheme. Practical issues to be aware of are: (i) that sharp changes in ridge direction will be difficult to resolve because the tangent calculation will be noisy and the transverse search interval of consecutive points on the ridge may intersect; and (ii) in situations where there is a complex system of closely aligned ridges, the ridge refinement may jump from one ridge to another if the search window overlaps nearby ridges. Typically, these issues can be resolved on a case-by-case basis by decreasing the step size in ridge tracking and refining over a smaller window many times.

\subsection{Results} 
\label{sec:ridge_results}

\begin{figure}
	\centering
		\includegraphics[width=6in]{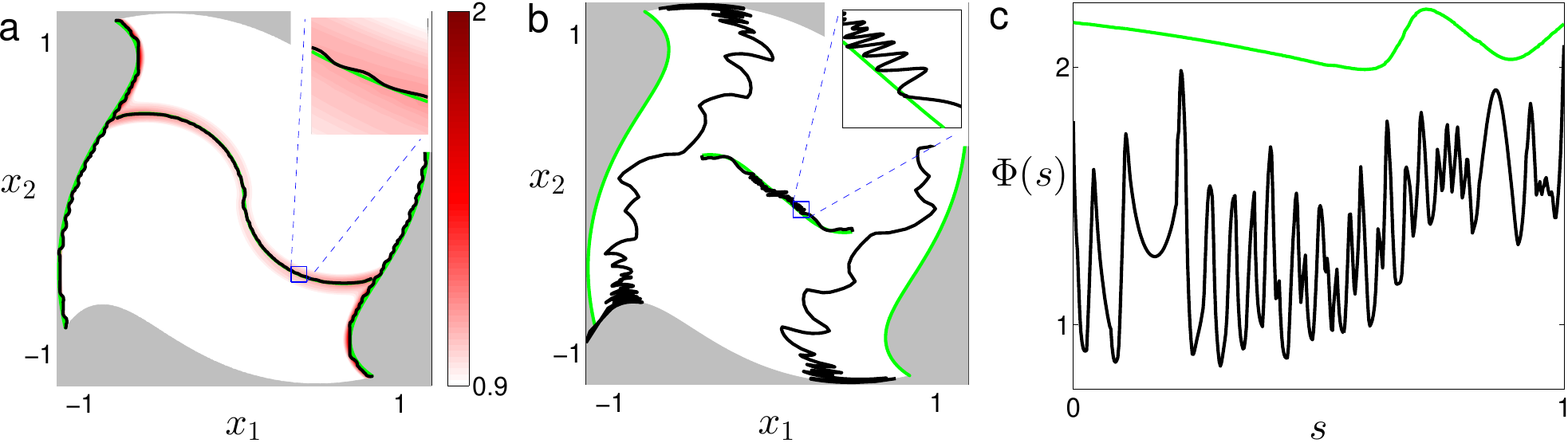}
		\caption{(a) Tracked ridges (black), and the corresponding refined ridges (green) of the FTLE field (red) for the autonomous system.  An inset is included to demonstrate the difference between the results. (b) Results of advection for the tracked and refined ridges.  The inset highlights how the small initial differences grow for the central ridge. (c) FTLE values, $\Phi(s)$, along the vertical tracked and refined ridge as a function of the in-line coordinate, $s$. }
	\label{fig:tracking_results}
\end{figure}

To demonstrate the ridge tracking and refinement approach, we apply these techniques to our model system.  We used a discretized velocity field with $\Delta x=2^{-11}$, relative tolerance of $10^{-7}$ for {\tt ode45}, and cluster size $\delta a=10^{-6}$.  The FTLE field is calculated on a uniform field of resolution $\delta x=0.01$.  The ridges found using the ridge tracking method are presented as black lines in figure~\ref{fig:tracking_results}(a). These do not smoothly follow the FTLE field due to the underlying discretization; this is clearly seen in the ridges near the left and right boundaries and on closer inspection of the center ridge.  Implementing the refinement scheme eliminates these artificial features, demonstrated by the green refined ridges that do smoothly follow the FTLE field.

To demonstrate the importance of refining the FTLE ridge, in figure~\ref{fig:tracking_results}(b) we present results for the advection of the unrefined and refined ridges.  The unrefined ridges along the left and right curved boundaries are actually advected into the domain, rather than remaining on the boundary as they should, and the unrefined ridge across the central region of the domain becomes jagged. In contrast, the refined ridges behave as expected and remain smooth. The second benefit of refinement is presented in figure~\ref{fig:tracking_results}(c).  Here the FTLE value is presented for the left ridge as a function of its in-line coordinate, $s$.  There are large oscillations in the FTLE value and while the result does appear to be continuous there are no features in the system of this scale, so it is clear that these variations are simply a result of an insufficiently refined ridge.  Refinement not only smooths out the results significantly, but we see that the FTLE ridge values are uniformly higher.  Similar results were also obtained for the central ridge.  

\section{FTLE ridge classification}
\label{sec:classification}

Because the FTLE field only represents the magnitude of the largest eigenvalue of the CG tensor, there is ambiguity regarding the deformation associated with an FTLE ridge, which, if available, would reveal a great deal about its fate and influence. A first attempt at FTLE classification defined three metrics measuring the normal repulsion, shear, and tangential stretching of the ridges~\cite{Tang2011}; while these metrics provide some of the desired information, the implementation relied on mapping the Hessian of the FTLE forward in time, which is numerically challenging.  In this section, we therefore further these ideas and modify them for application to refined FTLE ridges.  

\subsection{Methods}
\label{sec:class_method}

\begin{figure}
	\centering
		\includegraphics[width=5in]{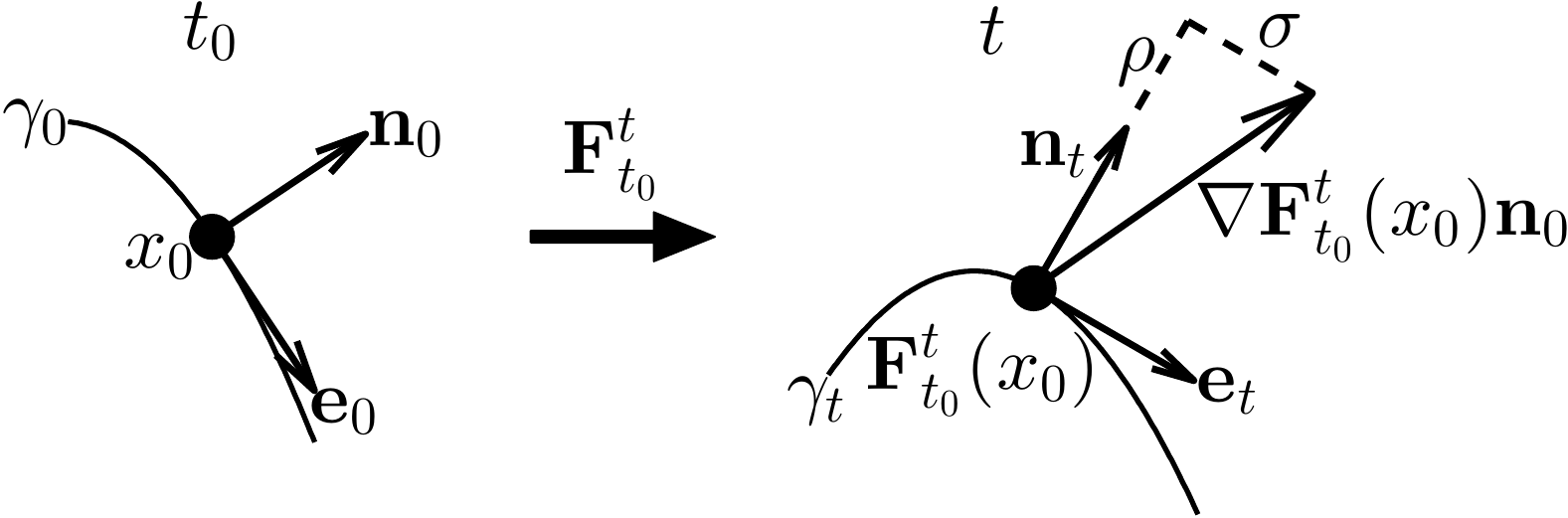}
		\caption{Sketch of the advection of a material line ($\gamma_0$) and the normal ($\bm{n}_0$) and tangent ($\bm{e}_0$) vectors at a point ($\bm{x}_0$) on the line.}
	\label{fig:adv_sketch}
\end{figure}

Given a material line, $\gamma_0$, at initial time, $t_0$, for each point along the line, $\bm{x}_0 \in \gamma_0$, the tangent and the normal vectors are $\bm{e}_0$ and $\bm{n}_0$, respectively.  As presented in figure~\ref{fig:adv_sketch}, the material line and its infinitesimal vectors are mapped forward by the flow map.  By definition, the tangent vector remains tangent to the material line, $\bm{\nabla F}^{t}_{t_0} \bm{e}_0 \parallel \bm{e}_t$, though stretching or contraction of the unit vector is possible.  Deformation of the initial normal vector, however, is such that it need not remain normal, and so the advected normal can be decomposed into two parts:
\begin{equation}
\label{eq:basic_advected_normal}
\bm{\nabla F}^{t}_{t_0}\bm{n}_0 = \rho \bm{n}_t + \sigma \bm{e}_t,
\end{equation}
where $\bm{n}_t$ and $\bm{e}_t$ are the normal and tangent vectors of the advected material line, respectively, $\rho$ is the normal stretching~\cite{Haller2011} and $\sigma$ is the Lagrangian shear~\cite{Haller2012}.

To classify the influence of the FTLE ridges, we consider the relative strengths of normal versus tangential growth, and the relative strengths of normal stretching and Lagrangian shear.  In this manner, through mapping the vectors $\bm{e}_0$ and $\bm{n}_0$, we first identify if material initially along or normal to the ridge is more greatly influenced by the local deformation.  Formally, we can quantify the magnitudes of the advected normal and tangential vectors via,
\begin{align}
\label{eq:ng}
n_l(\bm{x}_0) &= \log \big|\big| \bm{\nabla F}^{t}_{t_0}(\bm{x}_0)\bm{n}_0(\bm{x}_0) \big|\big|, \\
\label{eq:eg}
e_l(\bm{x}_0) &= \log \big|\big| \bm{\nabla F}^{t}_{t_0}(\bm{x}_0)\bm{e}_0(\bm{x}_0) \big|\big|,
\end{align}
where we have added the subscript $l$ to denote the logarithm scaling.  Because the logarithm scaling is used, stretching and contraction of the initial unit vectors are represented by positive and negative values, respectively. These metrics are related to the previously mentioned classification scheme\cite{Tang2011}, with the difference being that in this case both measures are concerned only with total growth and not their growth relative to the tangent and normal directions. In a similar manner, the hyperbolic repulsion and Lagrangian shear are given by:
\begin{align}
\label{eq:rl}
\rho_l &= \log \Big| \langle \bm{n}_t, \bm{\nabla F}^{t}_{t_0} \bm{n}_0 \rangle \Big| \\
\label{eq:sl}
\sigma_l &= \log \Big| \langle \bm{e}_t, \bm{\nabla F}^{t}_{t_0} \bm{n}_0 \rangle \Big|.
\end{align}
These metrics do match the previous classifications, except for a time averaging, the difference being in how $\bm{e}_t$ and $\bm{n}_t$ are calculated.  Note that the positive and negative values of $\rho_l$ correspond to repulsion and attraction in the normal direction, and that the convention (i.e. clockwise or counter-clockwise) of shear is lost by the absolute values and only the magnitude of shear is measured.  The key to calculating these quantities is to accurately identify the ridges so that advection is possible, enabling $\bm{e}_t$ to be reliably calculated by applying the flow map gradient to $\bm{e}_0$; $\bm{n}_t$ is then found by virtue of it being perpendicular to $\bm{e}_t$.  

As we have advocated throughout the paper, in performing any such calculations it is important to understand the inherent limitations given that FTLE ridges are the most unstable transport features in the flow field. Because the eigenvectors of the flow map gradient, which we refer to as $\bm{\xi}_1$ and $\bm{\xi}_2$, with corresponding eigenvalues $\lambda_2>\lambda_1>0$, are orthogonal, they provide a natural basis to decompose the normal and tangential ridge vectors and assess the impact of small errors.  As demonstrated in appendix~\ref{app:sensitivity}, when the normal and tangential ridge vectors are not closely aligned with $\bm{\xi}_1$ and $\bm{\xi}_2$, small errors do not significantly impact the classification metrics.  If the tangent and normal vector are approximately aligned with the eigenvector fields, however, small errors become a concern.  If an FTLE ridge is effectively a strainline (i.e. everywhere tangent with $\bm{\xi}_1$) errors do not greatly impact $n_l$ but can impact $e_l$, though only large amounts of error can cause the relative size of these two to be significantly altered. Furthermore, for FTLE ridges that are close to being strainlines, both $\rho_l$ and $\sigma_l$ are sensitive to errors, particularly when $\lambda_2 \gg \lambda_1$.  Full derivations of the sensitivity of the classification metrics that underly this summary are presented in appendix~\ref{app:sensitivity}.  This sensitivity analysis is particularly enlightening when considering ridges that closely align with the eigenvector fields where the classification metrics are most sensitive.

\subsection{Results}

\begin{figure}
	\centering
		\includegraphics[width=5in]{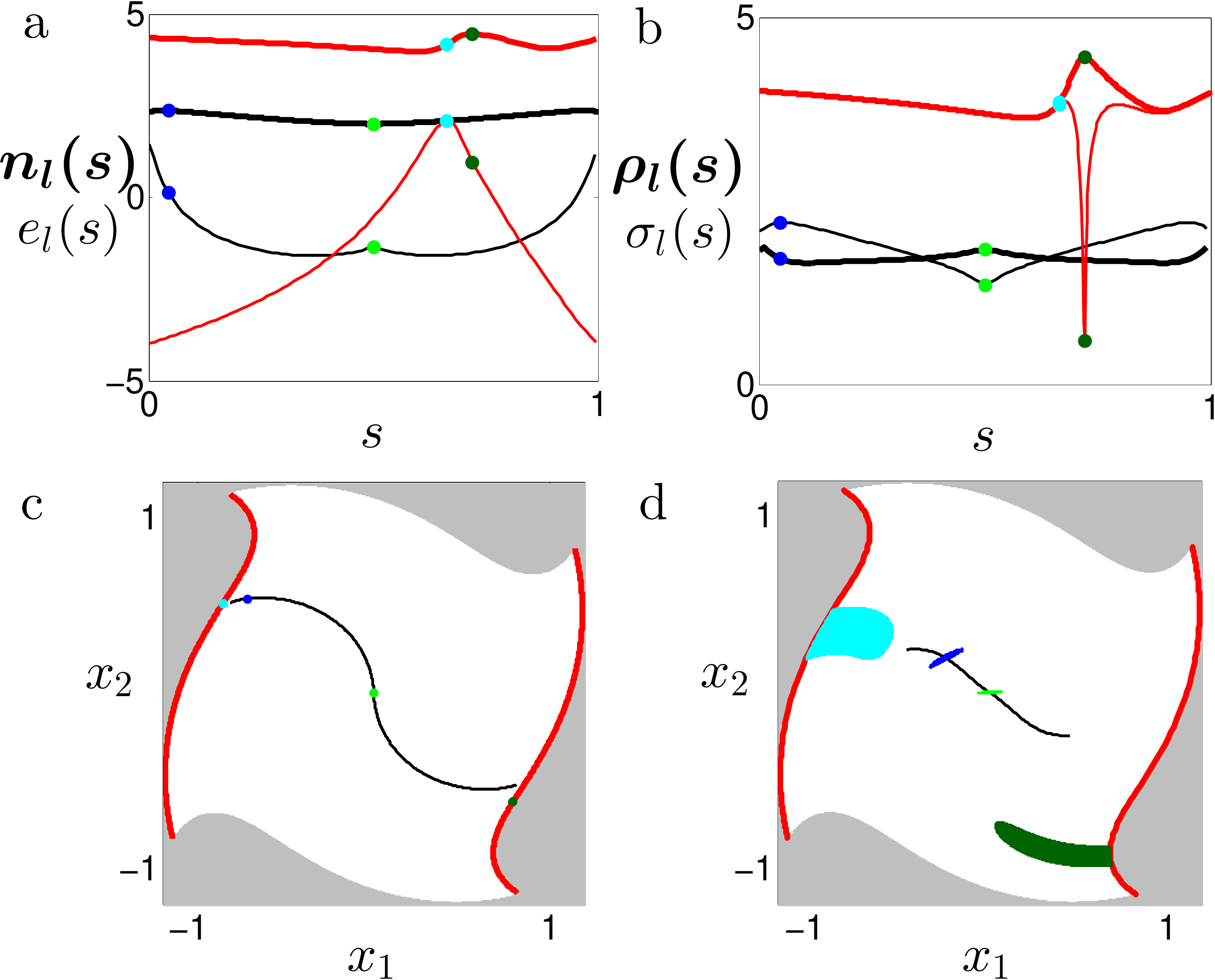}
		\caption{(a) Normal (thick lines) and tangential (thin lines) growth along the left (which is symmetric to the right) boundary (red liness) and center (black lines) FTLE ridges as a function of a normalized inline coordinate, s. (b) Normal repulsion (thick) and Lagrangian shear (thin) along the FTLE ridges. Advection of FTLE ridges and passive tracer patches from time (c) $t=0$ to (d) $t=2$.}
	\label{fig:classification}
\end{figure}

The classification metrics \eqref{eq:ng}, \eqref{eq:eg}, \eqref{eq:rl} and \eqref{eq:sl} are applied to the analytical model system, using the same resolution of discretized data used for the results in figure~\ref{fig:tracking_results}.  The classification values as a function of a normalized line coordinate, $s$, are presented in figures~\ref{fig:classification}(a) and~\ref{fig:classification}(b).  For this system, the values of $\lambda_1$ and $\lambda_2$, and the orientation of the ridges relative to the eigenvector field, all lend themselves to reliable classification.  The ridges along the boundary have greater normal than tangential growth along their entire length; indeed, only where these boundary ridges are near to the central ridge does their tangential growth become positive, and the negative values elsewhere indicate tangential contraction.  Similarly, there is growth of the normal vector everywhere along the ridge across the center of the domain, and tangential contraction everywhere except near the left and right boundary. For all three ridges, the relative amounts of normal repulsion and Lagrangian shear are approximately equal, with the exception that near the maximum value of $n_l$ along the boundary ridges, normal repulsion dominates. 

With this information in hand, we now know that as the central ridge is advected there will be contraction of material initially distributed along the ridge, stretching of material initially normal to the ridge, and that this latter stretching will be evenly distributed between shear and normal repulsion. Thus, we expect a circular patch of particles initially placed anywhere along the ridge to become a tilted ellipse with its semi-major axis oriented at some clearly observable angle to the ridge. This behavior is illustrated by the light green and dark blue patches in figure figures~\ref{fig:classification}(c) and~\ref{fig:classification}(d).  For the ridges along the boundaries, there is a segment where $n_l,e_l>0$ indicating that there will be stretching both along the tangent and normal to the tangent. As demonstrated by the dark green patches in figures~\ref{fig:classification}(c) and~\ref{fig:classification}(d), we see there is a region of very low Lagrangian shear so that the patch becomes stretched normally to the ridge without any significant tilting. In contrast, the light blue patch of particles released in the vicinity of the intersections between the boundary and central ridges experiences both tangential and normal growth; the growth of both components is reflected in the non-orthogonal intersection of the patch and ridge and the growth of the patch away from the ridge. The only condition not presented in this model is where $e_l>n_l$, which would result in an initially circular patch stretching along a ridge without having significant expansion of material initially normal to the ridge.

\section{Application to Ningaloo peninsula}
\label{sec:ningaloo}

As a further application of our FTLE-based methods, we investigated a data set produced by a numerical ocean model of a region containing the Ningaloo peninsula in Western Australia, which is the location of one of the world's longest fringing coral reefs and vast offshore hydrocarbon reserves.   The sea surface velocity fields employed in the study were produced by a double nested validated ROMS model with 2000m resolution, and we consider a 108-hour time window. A previous LCS analysis of this region investigated the impact of incorporating surface wind effects into strainline LCS analysis~\cite{Allshouse2015}.  We now calculate, advect, and classify the FTLE ridges of this system.  

Before analyzing the FTLE field, we performed a sequence of convergence studies to confirm reliable results.  First the relative tolerance of the ODE solver was varied to ensure that the flow map was calculated accurately, and the relative tolerance of $10^{-7}$ was sufficient for converged trajectory calculations.  Next, cluster sizes were varied to find converged values of the FTLE field.  Cluster sizes of 1m, $10^{-2}$m, and $10^{-4}$m were tested and the results for all three cluster sizes were nearly identical, so a cluster size of $1$m was used. While this may seem small given a grid resolution of 2000m, it should be noted that the ratio of grid resolution to cluster size is approximately the same as that used for the discretized studies of our analytical model.  As a final test, we subsample the domain to ensure that the results are not significantly affected. When decreasing the domain resolution to 4000m, we found that the dominant features of our analysis were nearly identical, with only slight (O(5\%)) differences in absolute value of the FTLE maxima.  Having determined the parameters needed for a robust FTLE calculation, the forwards- and backwards-time FTLE fields on a 200m grid were calculated in order to perform the initial ridge extractions, followed by ridge refinement and classification. 

\begin{figure}
	\centering
		\includegraphics[width=6.5in]{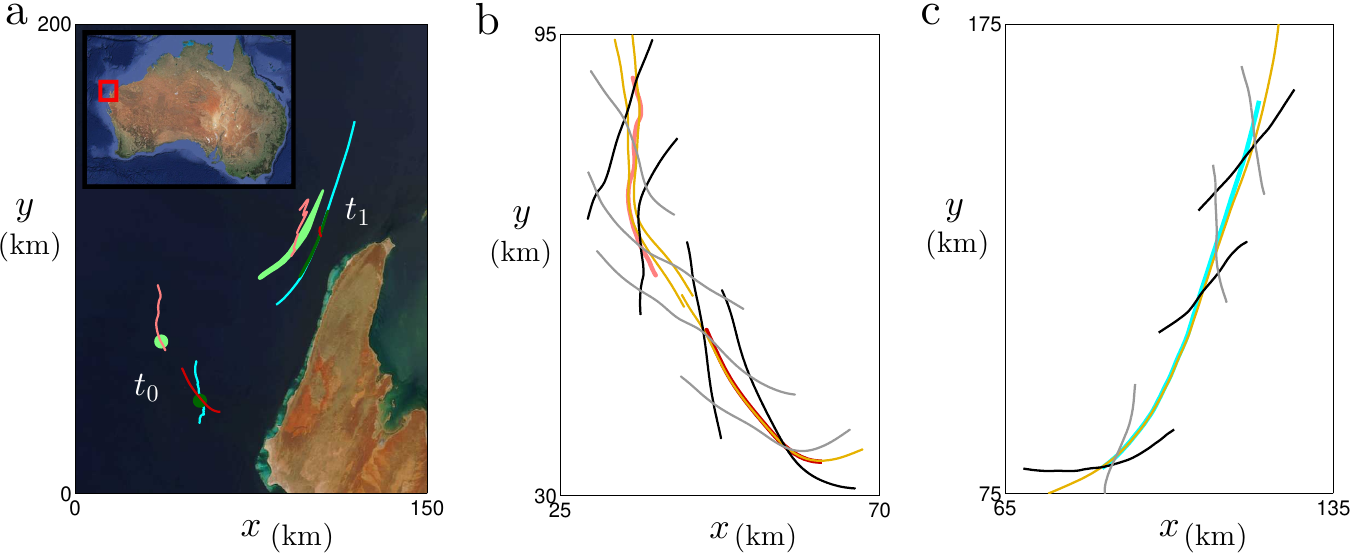}
		\caption{(a) Forward-time (light and dark red) and the backward FTLE ridges (light blue) along with passive tracer patches (light and dark green) at times $t_0=$0:00 21 December 2009 and $t_1=$ 12:00 25 December 2009 overlaying a MODIS satellite image of the Ningaloo peninsula.  Inset of Australia shows the location of the region.  Positive and negative shearlines (black and gray), strainlines (yellow), and (b) forward-time FTLE ridges or (c) backward-time FTLE ridges.}
	\label{fig:ningaloo_ridges}
\end{figure}

Two forward-time ridges that are in close proximity and a single backward time ridge were particularly strong features of the FTLE field, and so we focus our attention on these. The positions of all three ridges at the beginning and end of the time window, which required advection of the refined ridges, are presented in figure~\ref{fig:ningaloo_ridges}(a).  Particle patches have been added to help demonstrate the local deformations near the ridges. A Lagrangian hyperbolic point lies at the intersection of one of the forward-time FTLE ridges and the backward-time FTLE ridge. The other forward time FTLE ridge, however, was not associated with any backward time ridge. Figures~\ref{fig:ningaloo_ridges}(b) and (c) plot strainlines and shearlines in the vicinity of these FTLE ridges, and we see that the two FTLE ridges that form the Lagrangian hyperbolic point are nicely captured by strainlines; the isolated forward-time FTLE ridge, however, is not well captured by a single strainline or shearline. 

\begin{figure}
	\centering
		\includegraphics[width=6.5in]{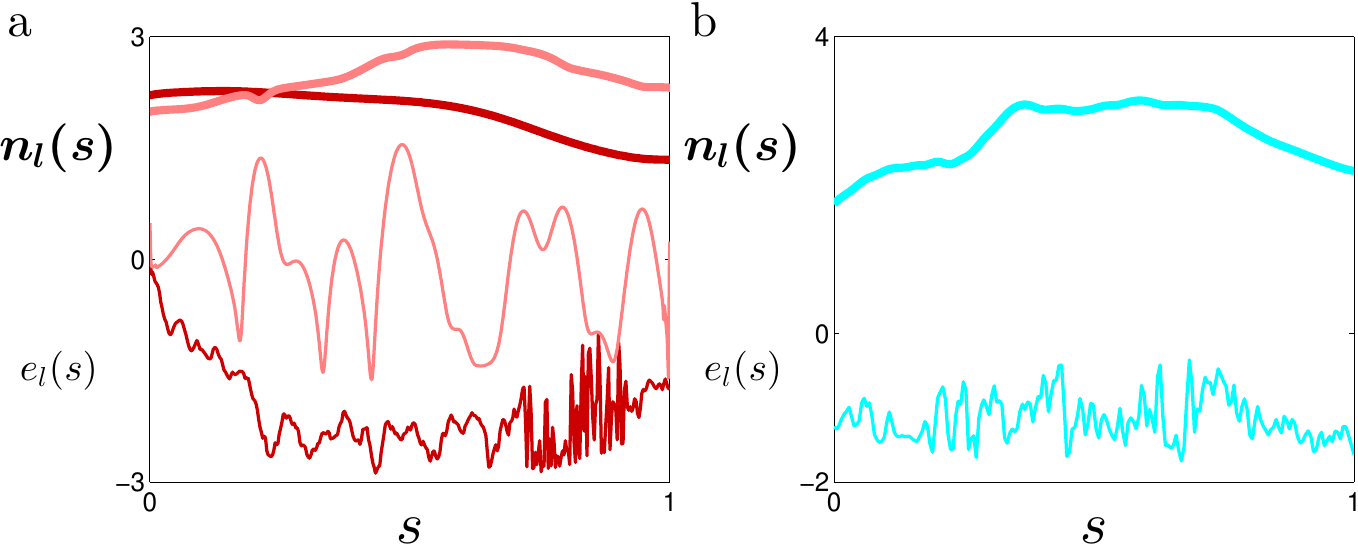}
		\caption{(a) $n_l$ (bold) and $e_l$ (thin) as a function of $s$ for the isolated forward-time ridge (light red) and intersecting forward-time ridge (dark red).  (b) $n_l$ and $e_l$ along the backwards-time ridge.}
\label{fig:ningaloo_classification}
\end{figure}

The classification scheme is implemented on these three ridges, and the values of $e_l$ and $n_l$ are presented as a function of a dimensionless inline coordinate $s$ in figure~\ref{fig:ningaloo_classification}(a) and \ref{fig:ningaloo_classification}(b) ($s=0$ corresponds to the right most end of the ridge).  For all three ridges there is growth of the normal vector along the entire length of the ridge, indicating that the ridge is influencing fluid on either side.  The persistent negative values of $e_l$ along the backward-time ridge means that this ridge will stretch in tangential length from time $t_0$ to $t_1$ (contraction in backwards time).  The $e_l$ value of the forward-time ridge that intersects the backwards time ridge is also persistently negative indicating that it will contract towards the backward-time ridge.  Interestingly, the third FTLE ridge has regions where $e_l<0$ and $e_l>0$, indicating there are alternating segments of expansion and contraction.  Additionally, we note that the values of $e_l$ along the two intersecting ridges appears to be noisy while the isolated forward-time ridge has a smoothly varying $e_l$.  This is indicative of the fact that the two intersecting ridges are well represented by strainlines.  

\begin{figure}
	\centering
		\includegraphics[width=6.5in]{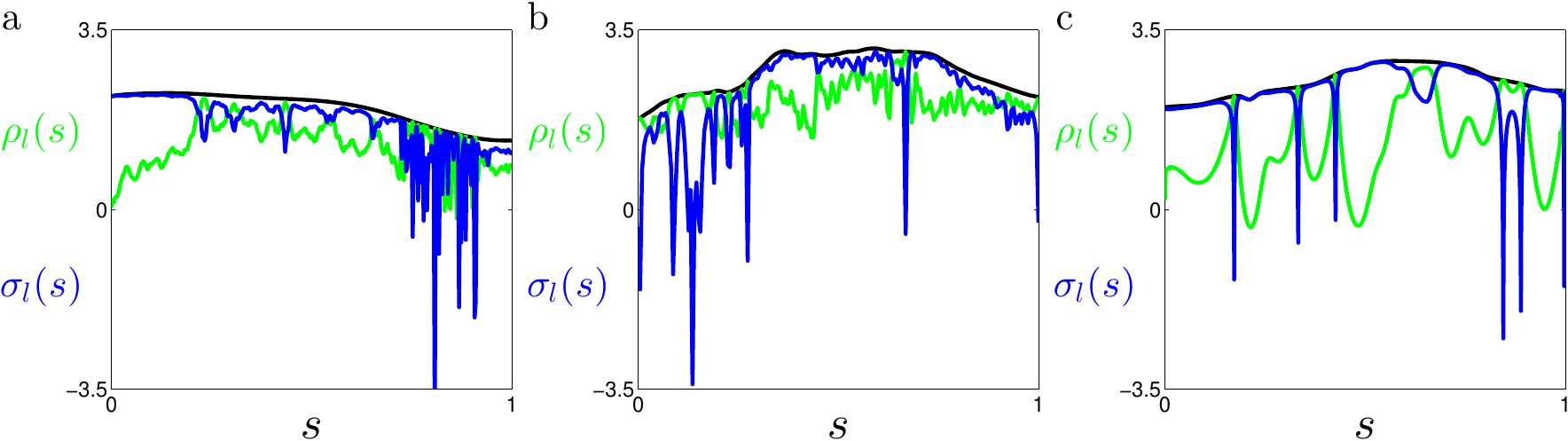}
		\caption{$\rho_l$ (green), $\sigma_l$ (blue), and $n_l$ (black) as a function of the scaled inline coordinate, $s$ for the (a) intersecting forward-time ridge, (b) intersecting backward-time ridge, and the (c) isolated forward-time ridge.}
\label{fig:ningaloo_rho}
\end{figure}

Finally, the values of $\rho_l$ and $\sigma_l$ are investigated along each of the three ridges with the results presented in figure~\ref{fig:ningaloo_rho}.  It should be noted that in all three cases the value of $n_l$ is smoothly varying along the entire length of the line (presented in black for reference).  The intersecting forward and backward-time ridges, presented in~\ref{fig:ningaloo_rho}(a) and~\ref{fig:ningaloo_rho}(b), respectively demonstrate large variation in the values of $\rho_l$ and $\sigma_l$ along the line.  This error is a manifestation of the sensitivity of these classification metrics when a ridge closely aligns with the CG eigenvector field (see Appendix B for further details).  The smoothly varying $\rho_l$ and $\sigma_l$ of the isolated ridge are presented in figure~\ref{fig:ningaloo_rho}(c).  While there are a couple isolated points of small Lagrangian shear, most of the ridge features larger values of Lagrangian shear than hyperbolic repulsion.  This manifests itself in the tilting of the particle patch relative to the ridge in figure~\ref{fig:ningaloo_ridges}.  These results are more reliable than the intersecting ridge classifications because in this instance, the ridges are not closely aligned with the eigenvector fields placing it in a regime where the metrics are robust to small errors.

\section{Conclusions}
\label{sec:conclusions}

While the FLTE method has its known limitations, it nevertheless identifies a particular type of LCS, namely the most active material lines in a fluid flow, and is a useful tool for preliminary analysis of flow transport. Of particular importance, therefore, is that the limitations of what can be interpreted from the FTLE field are clear and that proper care is taken to ensure that the FTLE field and its ridges are accurately determined. 

In this paper, we used both the finite-difference and advected-gradient approaches for calculating the FTLE field, and found that any benefits of the advected-gradient method are lost when working with discretized data. While it is reassuring that the finite-difference method is found to be robust, convergence of the flow map by varying the relative tolerance and of the finite-differences by varying cluster size are necessary. The normal-maximum ridge definition was introduced as a combination of the strengths of both the height and watershed ridges.  Recognizing the ridge detection scheme by~\citet{Lipinski2010}, minor modifications were made to take advantage of the parallelization of Matlab; a practical scheme was introduced to enable rapid refinement of an FTLE ridge.

The refined FTLE ridges mark material lines with large associated deformation, but classification and advection is furthermore necessary to understand how these material lines are influencing flow transport.  Measuring the growth of the normal and tangent vectors determines if an FTLE ridge is primarily influencing material along its length or adjacent to it.   In cases where normal growth dominates, it is then reasonable to compare the normal repulsion and Lagrangian shear.  We performed some analysis to determine when these different types of deformation are  sensitive to errors in the FTLE ridge. As a reasonable test, these techniques were applied to both an analytic model and an ocean surface data set.

\begin{acknowledgments}
The authors would like to thank Jean-Luc Thiffeault for discussions about the advected-gradient method, Doug Lipinski about ridge tracking, and George Haller about ridge definitions. This work was funded by ONR Grant Number N000141210665.

\end{acknowledgments}

\appendix

\section{Analytic model}
\label{app:model}

The following is a basic analysis of the autonomous system in order to get the FTLE field.  In order to calculate the FTLE field, it is necessary to obtain the flow map and its gradient in terms of the transformed coordinates.  First we can solve equations~\eqref{eq:base_equations} by separation of variables yielding:
\begin{align}
\label{eq:basic_fm}
X_1(t,A_1) &= \mbox{sign}(A_1) \left[ 1 - \left( 1 - \frac{1}{A_1{}^2} \right) e^{-2t} \right]^{-1/2}, \nonumber \\
X_2(t,A_2) &= \mbox{sign}(A_2) \left[ 1 - \left( 1 - \frac{1}{A_2{}^2} \right) e^{ 2t} \right]^{-1/2}.
\end{align}
where $(A_1,A_2)$ are the Lagrangian coordinates, i.e. initial positions, of the trajectories.  Then to get the flow map of the transformed system, we use substitution of the~\eqref{eq:basic_fm} into~
\eqref{eq:transform_cord}.  Because it is convenient to present the flow map in terms of its initial conditions in the transformed coordinates, we note that 
\begin{align}
A_1 &= a_1 \cos(r_0) + a_2\sin(r_0), \nonumber \\
A_2 &= a_2 \cos(r_0) - a_1\sin(r_0), 
\end{align}
where $(a_1,a_2)$ are the Lagrangian coordiantes in the transformed system and $r_0=\sqrt{A_1{}^2+A_2{}^2}=\sqrt{a_1{}^2+a_2{}^2}$.

Next to get the flow map gradient, $\bm{\nabla F}^t_{t_0}=\partial x_i/\partial a_j$, the chain rule is applied to ~\eqref{eq:transform_cord}.  Taking the derivative of equation~\eqref{eq:transform_cord} in terms of the Lagrangian coordinate $a_i$ gives:
\begin{align}
\frac{\partial x_1}{\partial a_i} &= \frac{\partial}{\partial a_i} \left[X_1\cos(r) - X_2\sin(r) \right], \nonumber \\
\label{appeq:der1}
&= \left[\frac{\partial X_1}{\partial a_i} - X_2 \frac{\partial r}{\partial a_i} \right] \cos(r) - \left[\frac{\partial X_2}{\partial a_i} + X_1\frac{\partial r}{\partial a_i}\right]\sin(r), \\
\frac{\partial x_2}{\partial a_i} &= \frac{\partial}{\partial a_i}\left[ X_2\cos(r) + X_1\sin(r)\right], \nonumber \\
\label{appeq:der2}
&= \left[\frac{\partial X_2}{\partial a_i} + X_1\frac{\partial r}{\partial a_i} \right] \cos(r) + \left[\frac{\partial X_1}{\partial a_i} - X_2 \frac{\partial r}{\partial a_i} \right]\sin(r).
\end{align}
Equations \eqref{appeq:der1} and \eqref{appeq:der2} indicate that we need to start deriving the terms $\partial X_i / \partial a_j$ and $\partial r / \partial a_i$.  From solving equation \eqref{eq:transform_cord} with the given initial conditions, it can be shown that
\begin{equation}
\label{appeq:dXdA}
\frac{\partial X_i}{\partial A_j} = \left[ \begin{array}{c c}
\frac{e^{-2t}}{A_1{}^3} \left[ 1 - \left( 1 - A_1{}^{-2} \right) e^{-2t}\right] ^{-3/2} & 0 \\
0 & \frac{e^{2t}}{A_2{}^3} \left[ 1 - \left( 1 - A_2{}^{-2} \right) e^{2t}\right]^{-3/2}
\end{array} \right],
\end{equation}
which can be utilized in the chain rule expansion to show that:
\begin{equation}
\label{appeq:dXda}
\frac{\partial X_i}{\partial a_j} = \frac{\partial A_k}{\partial a_j} \frac{\partial X_i}{\partial A_k}.
\end{equation}
If we reconfigure the coordinate transform it can be shown that
\begin{align*}
A_1 &= a_1\cos(r_0) - a_2\sin(r_0), \nonumber \\
A_2 &= a_2\cos(r_0) + a_1\sin(r_0), 
\end{align*}
where $r_0 = \sqrt{a_1{}^2+a_2{}^2}$.  This leads to
\begin{equation}
\label{appeq:dAda}
\frac{\partial A_i}{\partial a_j} = \left[ \begin{array}{c c}
\left(1-\frac{a_1a_2}{r_0}\right)\cos(r)-\frac{a_1{}^2}{r_0}\sin(r_0) & -\left(1+\frac{a_1a_2}{r_0}\right)\sin(r)-\frac{a_2{}^2}{r_0}\cos(r_0) \\
\left(1-\frac{a_1a_2}{r_0}\right)\sin(r)+\frac{a_1{}^2}{r_0}\cos(r_0) & \left(1+\frac{a_1a_2}{r_0}\right)\cos(r)-\frac{a_2{}^2}{r_0}\sin(r_0)
\end{array} \right],
\end{equation}
which when substituted along with equation~\eqref{appeq:dXdA} give us the necessary terms to calculate equation~\eqref{appeq:dXda}.

Calculating the remaining derivatives of equations~\eqref{appeq:der1} and \eqref{appeq:der2} is via
\begin{equation}
\label{appeq:drda}
\frac{\partial r}{\partial a_i} = \frac{\partial X_j}{\partial a_i}\frac{\partial r}{\partial X_j} = \frac{\partial X_j}{\partial a_i}\frac{X_j}{r}.
\end{equation}
where
\begin{displaymath}
\frac{\partial r}{\partial X_i} = \frac{X_i}{r}.
\end{displaymath}

With the necessary derivatives to calculate the terms of the flow map gradient, we can now calculate the terms of the CG tensor:
\begin{equation}
\bm{C}^{t}_{t_0}(\bm{x}_0) = \left[ \begin{array}{c c}
C_{11} & C_{12} \\
C_{12} & C_{22}
\end{array} \right] = 
\left[ \begin{array}{c c}
\frac{\partial x_1}{\partial a_1}^2 + \frac{\partial x_2}{\partial a_1}^2 & \frac{\partial x_1}{\partial a_1} \frac{\partial x_1}{\partial a_2}+ \frac{\partial x_2}{\partial a_1} \frac{\partial x_2}{\partial a_2}  \\
\frac{\partial x_1}{\partial a_1} \frac{\partial x_1}{\partial a_2}+ \frac{\partial x_2}{\partial a_1} \frac{\partial x_2}{\partial a_2} & \frac{\partial x_1}{\partial a_2}^2 + \frac{\partial x_2}{\partial a_2}^2
\end{array} \right].
\end{equation}
Finally, with the generic terms of the CG tensor, the largest eigenvalue and eigenvector can be shown to be:
\begin{equation}
\label{appeq:swirl_analytic_l2}
\lambda_2 = \frac{C_{11} + C_{22} + \sqrt{ (C_{11} - C_{22})^2 + 4 C_{12}{}^2 }}{2}.
\end{equation}

\section{Sensitivity analysis}
\label{app:sensitivity}

Central to the following arguments is the singular value decomposition of the flow map gradient.  We define the smaller singular value as $\lambda_1$ and it has the corresponding right-singular vector, $\bm{\xi}_1$, which maps forward to the left-singular vector, $\bm{u}_1$ such that
\begin{align*}
\bm{\nabla F}^{t}_{t_0} \bm{\xi}_1 &= \lambda_1 \bm{u}_1, \\
\bm{\nabla F}^{t}_{t_0} \bm{\xi}_2 &= \lambda_2 \bm{u}_2. 
\end{align*}
Throughout the derivation we utilize the conventions $\bm{\xi}_1 \times \bm{\xi}_2 = \bm{u}_1 \times \bm{u}_2 = 1$.  

A point, $\bm{x}_0$, along the FTLE ridge, $\gamma_0$, is allowed to have the general orientation:
\begin{align*}
\bm{e}_0 &= \sqrt{1-b^2}\bm{\xi}_1 + b\bm{\xi}_2\\
\bm{n}_0 &= -b\bm{\xi}_1 + \sqrt{1-b^2}\bm{\xi}_2.
\end{align*}
where $|b|\leq1$ and the tangent and normal vectors have unit length.  In the unique case where the FTLE ridge is a strainline (stretchline), $b=0$ ($b=1$).  For the following analysis, a perturbation $\epsilon \ll 1$ is added to the tangent vector such that the $\bm{\xi}_2$ coefficient becomes $b+\epsilon$, and so the other components of this vector and of the normal vector are modified to  that both remain of unit length to leading order. This leads to the perturbed tangent and normal vectors:
\begin{align*}
\bm{e}_0' &= \left(\sqrt{1-b^2} -\epsilon\frac{b}{\sqrt{1-b^2}} - \epsilon^2 \frac{1}{2(1-b^2)^{3/2}} + O(\epsilon^3) \right)\bm{\xi}_1 + (b+\epsilon)\bm{\xi}_2,\\
\bm{n}_0' &= -(b+\epsilon)\bm{\xi}_1 + \left(\sqrt{1-b^2} -\epsilon\frac{b}{\sqrt{1-b^2}} - \epsilon^2 \frac{1}{2(1-b^2)^{3/2}}+ O(\epsilon^3) \right)\bm{\xi}_2.
\end{align*}
where $'$ denotes the perturbed vector.  Mapping these forward under the action of the flow map gradient gives:
\begin{align*}
\bm{\nabla F}^{t}_{t_0}\bm{e}_0' &= \left(\sqrt{1-b^2} -\epsilon\frac{b}{\sqrt{1-b^2}} - \epsilon^2 \frac{1}{2(1-b^2)^{3/2}} + O(\epsilon^3) \right)\lambda_1\bm{u}_1 + (b+\epsilon)\lambda_2\bm{u}_2,\\
\bm{\nabla F}^{t}_{t_0}\bm{n}_0' &= -(b+\epsilon)\lambda_1\bm{u}_1 + \left(\sqrt{1-b^2} -\epsilon\frac{b}{\sqrt{1-b^2}} - \epsilon^2 \frac{1}{2(1-b^2)^{3/2}}+ O(\epsilon^3) \right)\lambda_2\bm{u}_2,
\end{align*}
and it can be shown that:
\begin{align*}
e_l' &= \lambda_2 \left[\delta^2 + (b + \epsilon)^2(1 - \delta^2) + \delta^2 O(\epsilon^3) \right]^{1/2}, \\
n_l' &= \lambda_2 \left[1        - (b + \epsilon)^2(1 - \delta^2) +  O(\epsilon^3) \right]^{1/2}.
\end{align*}
The advected tangent and normal are thus:
\begin{align*}
\bm{e}_t' &= \frac{1}{e_l'}\left[\left(\sqrt{1-b^2} -\epsilon\frac{b}{\sqrt{1-b^2}} - \epsilon^2 \frac{1}{2(1-b^2)^{3/2}} + O(\epsilon^3) \right)\lambda_1\bm{u}_1 + (b+\epsilon)\lambda_2\bm{u}_2\right]\\
\bm{n}_t' &= \frac{1}{e_l'}\left[-(b+\epsilon)\lambda_2\bm{u}_1 + \left(\sqrt{1-b^2} -\epsilon\frac{b}{\sqrt{1-b^2}} - \epsilon^2 \frac{1}{2(1-b^2)^{3/2}} + O(\epsilon^3) \right)\lambda_1\bm{u}_2\right].
\end{align*}
The normally hyperbolic repulsion and Lagrangian shear can be determined via:
\begin{align*}
\rho' &= \langle \bm{n}_t' , \bm{\nabla F}^{t}_{t_0}\bm{n}_0'\rangle = \frac{\lambda_1 \lambda_2}{e_l'} \left[1 + O(\epsilon^3)  \right], \\
\sigma' &= \langle \bm{e}_t' , \bm{\nabla F}^{t}_{t_0}\bm{n}_0' \rangle = \frac{\lambda_2^2- \lambda_1^2}{e_l'} \left[ b \sqrt{1-b^2} + \epsilon \frac{1-2b^2}{\sqrt{1-b^2}} - \epsilon^2 \frac{3b-2b^2}{2(1-b^2)^{3/2}} + O(\epsilon^3)  \right].
\end{align*}

To gain insight into the stability of these derivations to perturbation we need to expand the terms and look at the coefficient multiplying $\epsilon$.  Of particular interest are three cases.  The first case is strainlines ($b=0$); these correspond to material lines for which normal repulsion is maximized and  thus are tangent to the $\bm{\xi}_1$-field.  We also consider the stretchlines ($b=1$), for which stretching is principally in a tangential direction.  Finally, we consider cases where $\lambda_1 \ll \lambda_2$ corresponding to $\delta \ll 1$; this is a common scenario for systems with large amounts of stretching.

The first classification metric we look at is the growth of the tangent vector, $e_l$. The growth of the perturbed tangent vector is obtained from: 
\begin{align*} 
e_l' &= \lambda_2 [ \delta^2 + (b+\epsilon)^2(1-\delta^2)+\delta^2 O(\epsilon^3)]^{1/2} \\
&=  \lambda_2 [ \delta^2 + b^2(1 - \delta^2)]^{1/2} + \epsilon \frac{ b (1-\delta^2)}{\sqrt{\delta^2 + b^2(1-\delta^2)}}  + O(\epsilon^2)
\end{align*}
The coefficient for the perturbation remains small for all cases we consider.  It goes to $1-\delta^2$ as $b\rightarrow 1$; therefore, $e_l$ for a stretchline is not significantly altered when the tangent vector orientation is perturbed.  For $\delta \ll 1$, the perturbation coefficient goes to 1, so when the vector is not closely aligned to either eigenvector field, perturbations will not result in large changes in $e_l$. The perturbation coefficient becomes zero, however, when $b=0$, and so we expand to consider the second order perturbation.  It can be shown that when $b=0$, the second order expansion becomes:
\begin{displaymath}
e_l' = \lambda_1 \left[ 1 - \frac{1}{2} \epsilon^2 \left(\frac{1}{\delta^2} - 1 \right) \right].
\end{displaymath}
This shows that when $\delta \ll 1$ and $b=0$, the coefficient for the perturbation diverges unless $\epsilon \ll \delta$.  This is not surprising since, for a strainline, the tangent is aligned with the direction of smallest stretching, $\bm{\xi}_1$, and the normal is aligned with the direction of largest stretching, $\bm{\xi}_2$.  Any perturbation to these vectors adds a contribution from $\bm{\xi}_2$ to the tangent vector; this component grows significantly causing the overall length of the advected tangent vector to be much larger than the unperturbed state.

Next we consider the growth of a perturbed normal vector:  
\begin{align*} 
n_l' &= \lambda_2 [ 1 - (b+\epsilon)^2(1-\delta^2)+\delta^2 O(\epsilon^3)]^{1/2} \\
&=  \lambda_2 [1 - b^2(1- \delta^2)]^{1/2} + \epsilon \frac{ b (1-\delta^2)}{\sqrt{1 - b^2(1- \delta^2)}}  + O(\epsilon^2).
\end{align*}
For the stretchline case ($b=1$), the perturbation coefficient reduces to $1/\delta - \delta$.  When $\delta \ll 1$, this is very large indicating that $n_l$ for stretchlines is sensitive to small perturbations; this can be rationalized in much the same way that sensitivity of $e_l$ for strainlines was. The perturbation term again reduces to zero for the strainline case ($b=0$), so the expansion to second order is again performed:
\begin{displaymath}
n_l' = \lambda_2 \left[ 1 + \frac{1}{2}\epsilon^2 (1-\delta^2) \right].
\end{displaymath}
In this case, small perturbations do not significantly alter the $n_l$ value because while the absolute error may be large, relative to the unperturbed value it is not.  For values of $b$ not near the extremes, the perturbation coefficient is of order one meaning that the values are not sensitive to perturbation.

The sensitivity of the normal hyperbolic growth depends on the expansion of $1/e_l'$:
\begin{align*}
\rho_l' &= \frac{\lambda_1 \lambda_2}{e_l'} \\
&= \frac{\lambda_1 \lambda_2}{\lambda_2[\delta^2 + b^2 (1-\delta^2)]^{1/2}}\left[1 - \epsilon \frac{b (1-\delta^2)}{\delta^2 + b^2 (1-\delta^2)} \right]
\end{align*}
For $\delta \ll 1$ and $b\neq0$ the perturbation coefficient remains small.  In the case where $b=0$, we have shown that $e_l'$ is sensitive to perturbation so it is not surprising that $\rho_l'$ will also be sensitive to errors.  Through substitution, it can be shown that for $b=0$:
\begin{displaymath}
\rho_l' = \lambda_2 \left( 1 + \frac{1}{2} \epsilon^2 ( \frac{1}{\delta^2} - 1) \right).
\end{displaymath}
Unless $\epsilon \ll \delta$ this term is also sensitive to perturbations.

Finally, the analysis of the Lagrangian shear produces similar results to the hyperbolic repulsion.
\begin{align*}
\sigma_l' &= \frac{\lambda_2^2- \lambda_1^2}{e_l'} \left[ b \sqrt{1-b^2} + \epsilon \frac{1-2b^2}{\sqrt{1-b^2}} - \epsilon^2 \frac{3b-2b^2}{2(1-b^2)^{3/2}} + O(\epsilon^3)  \right] \\
&= \frac{(\lambda_2^2- \lambda_1^2)b\sqrt{1-b^2}}{\lambda_2[\delta^2 + b^2 (1-\delta^2)]^{1/2}} \left[ 1 + \epsilon \left(\frac{1-2b^2}{b(1-b^2)} - \frac{b (1-\delta^2)}{\delta^2 + b^2 (1-\delta^2)}\right) +O(\epsilon^2) \right]
\end{align*}
For $b=1$, there is a singularity in the perturbation coefficient indicating that perturbations for a stretchline will result in large changes in $\sigma_l'$.  Again for the strainline, it is necessary to expand to the second order to study the stability of perturbations.  This results in the expansion:
\begin{displaymath}
\sigma_l' = \frac{\lambda_2^2-\lambda_1^2}{\lambda_1}\left( \epsilon + \frac{1}{2}\epsilon^3\left(\frac{1}{\delta^2}-3\right)\right)
\end{displaymath}
It makes sense that this term is only insensitive to noise when $\epsilon \ll \delta$, because it relies on $e_l'$ being insensitive to noise.

\nocite{*}
\bibliography{AllshouseChaos}
\bibliographystyle{plainnat}

\end{document}